\documentclass[a4paper,english]{article}
\usepackage[T1]{fontenc}
\usepackage[latin9]{inputenc}
\usepackage{mathrsfs}
\usepackage{amsmath}
\usepackage{amssymb}



\usepackage{babel}
\begin{document}

\title{Distance, entropy and similarity measures of Type-2 soft sets}

\author{Rajashi Chatterjee, P. Majumdar, S. K. Samanta}
\maketitle
\begin{abstract}
The concept of Type-2 soft sets had been proposed as a generalization
of Molodstov's soft sets. In this paper some shortcomings of some
existing distance measures for Type-1 soft sets have been shown and
accordingly some new distance measures have been proposed. The axiomatic
definitions for distance, entropy and similarity measures for Type-2
soft sets have been introduced and a couple of such measures have
been defined. Also the applicability of one of the proposed similarity
measures have been demonstrated by showing its utility as an effective
tool in a decision making problem.
\end{abstract}

\section{Introduction }

Soft set theory, introduced by Molodstov \cite{molodstov}, emerged
as a revolutionary technique in dealing with intrinsic imprecision
with the help of adequate parameterization. Thereby the theory of
soft sets not only catered to the shortcoming encountered by fuzzy
sets \cite{zadehf1} possibly due to the lack of a parameterization
tool but in his work, Molodstov had also shown that fuzzy sets were
special types of soft sets, the parameter set being considered over
the unit interval $[0,1]$. Later, Maji et. al \cite{majisst} presented
a detailed theoretical study on soft sets. Thereafter the theory of
soft sets has undergone rapid developments in different directions
(\cite{uni}-\cite{pment},\cite{infsys}).\\
The organization of the paper is as follows:

Section 1 is the introductory portion, Section 2 is dedicated to recalling
some useful preliminary results. The notions of distance, entropy
and similarity measures for type-2 soft sets are defined in Sections
3, 4 and 5 respectively. In Section 6, an application of the proposed
similarity measure in a decision making problem is shown. Section
7 concludes the paper.

\section{Useful preliminaries}

Before introducing the proposed measures we first brush up some preliminary
results that would be useful for future purposes.\\
\\
\textbf{Definition} \textbf{2.1. }\cite{majisst} Let $X$ be an initial
universe and $E$ be a set of parameters. Let $\mathscr{P}(X)$ denote
the power set of $X$ and $A\subset E$. A pair $(F,A)$ is called
a soft set if $F$ is a mapping of $A$ into $\mathscr{P}(X)$.\\
\\
\textbf{Definition 2.2.} \cite{semiring} The bi-intersection of two
T1SS $(F,A)$ and $(G,B)$ is defined as $(F,A)\tilde{\cap}(G,B)=(H,C)$,
where $C=A\cap B$ such that for all $\alpha\epsilon C$, $H(\alpha)=F(\alpha)\cap G(\alpha)$.\\
\\
\textbf{Definition} \textbf{2.3.} $(F,A)\tilde{\cap}(G,B)=\begin{cases}
(H,A\cap B) & if\, A\cap B\neq\varphi\\
\widehat{\phi},\, the\, empty\, function & if\, A\cap B=\varphi
\end{cases}$ \\
where $H(\alpha)=F(\alpha)\cap G(\alpha)$, for all \textbf{$\alpha\epsilon A\cap B$.}\\
\textbf{}\\
\textbf{Definition} \textbf{2.4.}\cite{pment} A T1SS $(F,A)$ is
said to be deterministic if it satisfies the conditions\\
$(i)$ $\cup_{e\epsilon A}F(e)=X$\\
$(ii)$ $F(e)\cap F(f)=\varphi$, $e,f\epsilon A$, $e\neq f$.\\
\\
\textbf{Definition} \textbf{2.5. }\cite{real} Let $R$ be the set
of real numbers, $\mathscr{P}(R)$ be the collection of all non-empty
bounded subsets of $R$ and $A$ be a set of parameters. Then a mapping
$F:A\rightarrow\mathscr{P}(R)$ is called a soft real set, which is
denoted by $(F,A)$. In particular, if $(F,A)$ is a singleton soft
set, then identifying $(F,A)$ with the corresponding soft element,
it is called a soft real number.\\
\\
\textbf{Definition} \textbf{2.6.} \cite{me} Let $(X,E)$ be an initial
soft universe and $A,E_{A}\subset E$ be two sets of parameters. Suppose
that $S(X,E_{A})$ denotes the set of all T1SS over the soft universe
$(X,E_{A})$. Then a type-2 soft set over $(X,E)$ is a mapping $\mathscr{F}:A\rightarrow S(X,E_{A})$.
It is denoted by $[\mathscr{F},A]$.

In this case the set of parameters $A$ is referred to as the $primary\, set\, of\, parameters$
while the set of parameters $E_{A}$ is known as the $underlying\, set\, of\, parameters$.
\\
Thus, for a type-2 soft set, corresponding to each parameter $\alpha\epsilon A$,
there exists a TISS $(F_{\alpha},E_{A}^{\alpha})\epsilon S(X,E_{A})$
where $F_{\alpha}:E_{A}^{\alpha}\rightarrow\mathscr{\mathscr{P}}(X)$
and $E_{A}^{\alpha}\subset E_{A}$ such that $\cup_{\alpha\epsilon A}E_{A}^{\alpha}=E_{A}$.
\\
\\
\textbf{Example 2.7.} Let $X=\{h_{1},h_{2},h_{3},h_{4},h_{5}\}$ be
a set of five houses. Suppose\\
 $E=\{beautiful,single\, storeyed,wooden,in\, good\, repair,spacious,near\, the\,$\\
$market,in\, green\, surroundings,furnished,luxurious,with\, good\, security,$\\
$with\, pool\}$\\
Let $A=\{beautiful,luxurious\}$ and $E_{A}=\{wooden,in\, green\, surroundings,$\\
$with\, good\, security\}$ be two sets of parameters such that $A,E_{A}\subset E$.
A T2SS $[\mathscr{F},A]$ can be defined as follows:\\
\\
$\mathscr{F}(beautiful)=\{\frac{wooden}{\{h_{2},h_{5}\}},\frac{in\, green\, surroundings}{\{h_{1},h_{2},h_{3},h_{4}\}}\}$\\
$\mathscr{F}(luxurious)=\{\frac{wooden}{\{h_{5}\}},\frac{with\, good\, security}{\{h_{1},h_{3},h_{5}\}}\}$\\
\\
Another T2SS $[\mathscr{G},B]$ can be defined on the same soft universe
$(X,E)$ as,\\
\\
$\mathscr{G}(spacious)=\{\frac{wooden}{\{h_{5}\}},\frac{with\, pool}{\{h_{3},h_{5}\}}\}$\\
$\mathscr{G}(beautiful)=\{\frac{wooden}{\{h_{2},h_{5}\}},\frac{near\, the\, market}{\{h_{4}\}},\frac{in\, green\, surroundings}{\{h_{2},h_{3}\}}\}$\\
\\
In this case, $B=\{spacious,beautiful\}$ and\\
 $E_{B}=\{wooden,with\, pool,near\, the\, market,in\, green\, surroundings\}$
\\
\\
\textbf{Remark 2.8.\cite{me}} Type 2 fuzzy sets introduced by L.
A. Zadeh may be considered as special types of type-2 soft sets when
the parameters are considered over $[0,1]$. Let $X$ be the universe
under consideration. The point-valued representation of a type-2 fuzzy
set , denoted by $\tilde{A}$ , is given in \cite{mendel} as,\\
$\tilde{A}=\{((x,u),\mu_{\tilde{A}}(x,u))|\forall x\epsilon X,\forall u\epsilon J_{x}\subseteq[0,1]\}$\\
 where $\mu_{\tilde{A}}(x,u)$ is the type 2 membership function for
all $x\epsilon X$. Then, for any $\alpha\epsilon[0,1]$ , the corresponding
$\alpha$-plane of $\tilde{A}$ is given by,\\
$\tilde{A_{\alpha}}=\{(x,u),\mu_{\tilde{A}}(x,u)\geq\alpha|\forall x\epsilon X,\forall u\epsilon J_{x}\subseteq[0,1]$\}
\\
For a particular $\alpha\epsilon[0,1]$, construct the set,\\
$S_{\alpha}=\{u:(x,u)\epsilon\tilde{A_{\alpha}}\}$\\
 Define a mapping $F_{\alpha}:S_{\alpha}\rightarrow\mathscr{P}(X)$
by\\
$F_{\alpha}(u)=\{x:(x,u)\epsilon\tilde{A_{\alpha}}\},u\epsilon S_{\alpha}$
. \\
Thus, $(F_{\alpha},S_{\alpha})$ constitutes a type-1 soft set over
$(X,[0,1])$. \\
Next, define a mapping $\mathscr{F}:[0,1]\rightarrow S(X,[0,1])$
such that $\mathscr{F}(\alpha)=(F_{\alpha},S_{\alpha})$.\\
Thus, $[\mathscr{F},[0,1]]$ is a type-2 soft set over $(X,[0,1])$.
\\
\textbf{}\\
\textbf{Definition} \textbf{2.9.} The various operations over T2SS
are defined as,\\
\\
Containment: $[\mathscr{F},A]\sqsubseteq[\mathscr{G},B]$ iff $A\subseteq B$
and $\mathscr{F}(\alpha)\tilde{\subseteq}\mathscr{G}(\alpha)$, for
all $\alpha\epsilon A$.

Thus, when $[\mathscr{F},A]\sqsubseteq[\mathscr{G},B]$ we have, $\mathscr{F}(\alpha)\tilde{\subseteq}\mathscr{G}(\alpha)$,
for each $\alpha\epsilon A$ and so it follows that $E_{A}^{\alpha}\subseteq E_{B}^{\alpha}$
and $F_{\alpha}(\beta)=G_{\alpha}(\beta)$, for each $\alpha\epsilon A$
and $\beta\epsilon E_{A}^{\alpha}$. Thus,\\
$E_{A}=\cup_{\alpha\epsilon A}E_{A}^{\alpha}\subseteq\cup_{\alpha\epsilon B}E_{B}^{\alpha}=E_{B}$
\\
\\
Equality: $[\mathscr{F},A]=[\mathscr{G},B]$ iff $[\mathscr{F},A]\sqsubseteq[\mathscr{G},B]$
and $[\mathscr{G},B]\sqsubseteq[\mathscr{F},A]$.\\
\\
Union: $[\mathscr{F},A]\sqcup[\mathscr{G},B]=\begin{cases}
\mathscr{F}(\alpha) & \alpha\epsilon A-B\\
\mathscr{G}(\alpha) & \alpha\epsilon B-A\\
\mathscr{F}(\alpha)\tilde{\cup}\mathscr{G}(\alpha) & \alpha\epsilon A\cap B
\end{cases}$\\
$[\mathscr{F},A]\sqcup[\mathscr{G},B]$ is a T2SS with $A\cup B$
being the set of primary parameters and the set of underlying parameters
being $\cup_{\alpha\epsilon A\cup B}(E_{A}^{\alpha}\cup E_{B}^{\alpha})=E_{A\cup B}.$\\
\\
Intersection:\\
 $[\mathscr{F},A]\sqcap[\mathscr{G},B]=\begin{cases}
\mathscr{F}(\alpha)\tilde{\cap}\mathscr{G}(\alpha),\alpha\epsilon A\cap B & whenever\, A\cap B\neq\varphi\\
\widehat{\phi},\, the\, empty\, function & if\, A\cap B=\varphi,the\, null\, set
\end{cases}$.\\
Whenever $A\cap B\neq\varphi$, the intersection between two T2SS
is defined below:\\
\\
\textbf{Definition} \textbf{2.10.} A T2SS $[\mathscr{F},A]$, is said
to be an absolute T2SS if for each $\alpha\epsilon A$, $\mathscr{F}(\alpha)$
is an absolute T1SS. It is denoted by $\tilde{\mathbb{A}}$.\\
\\
\textbf{Definition} \textbf{2.11.} A T2SS $[\mathscr{F},A]$, is said
to be a null T2SS if for each $\alpha\epsilon A$, $\mathscr{F}(\alpha)$
is a null T1SS. It is denoted by $\tilde{\Phi}$.\\
 \\
\textbf{Definition} \textbf{2.12.} Two T2SS $[\mathscr{F},A]$ and
$[\mathscr{G},B]$ are said to be $disjoint$ if $A\cap B=\varphi$.
They are said to be $weakly\, disjoint$ if $A\cap B\neq\varphi$
but corresponding to each $\alpha\epsilon A\cap B$, $E_{A}^{\alpha}\cap E_{B}^{\alpha}=\varphi$
i.e. in such cases $[\mathscr{F},A]\sqcap[\mathscr{G},B]=\tilde{\varphi},$the
null T1SS. They are said to be $elementwise\, disjoint$ if $A\cap B\neq\varphi$,
$E_{A}\cap E_{B}\neq\varphi$ and for each $\alpha\epsilon A\cap B$,
$\mathscr{F}(\alpha)\tilde{\cap}\mathscr{G}(\alpha)$ is a null T1FS.
i.e. corresponding to each primary parameter $\alpha\epsilon A\cap B$,
$(F_{\alpha},E_{A}^{\alpha})\tilde{\cap}(G_{\alpha},E_{B}^{\alpha})=\tilde{\Phi}$
over the set of parameters $E_{A\cap B}$.

\section{Distance measures}

\textbf{Definition} \textbf{3.1. }\cite{kdist} Let $(F,A)$, $(G,B)$
and $(H,C)$ be any 3 soft sets in a soft space $(X,E)$ and $d:S(X,E)\times S(X,E)\rightarrow\mathbb{R}^{+}$,
a mapping. Then\\
$(1)$ $d$ is said to be a quasi-metric if it satisfies\\
$(M1)$ $d((F,A),(G,B))\geq0$\\
$(M2)$ $d((F,A),(G,B))=d((G,B),(F,A))$ \\
$(2)$ A quasi-metric $d$ is said to be a semi-metric if\\
$(M3)$$d((F,A),(G,B))+d((G,B),(H,C))\geq d((F,A),(H,C))$\\
$(3)$ A semi-metric $d$ is said to be a pseudo-metric if\\
$(M4)$ $(F,A)=(G,B)$ $\Rightarrow d((F,A),(G,B))=0$\\
$(4)$ A pseudo metric $d$ is said to be a metric if\\
$(M5)$ $d((F,A),(G,B))=0$ $\Rightarrow(F,A)=(G,B)$.\\
\\
Following this definition, Kharal\cite{kdist} had introduced several
new measures of calculating the distance between T1SS. In his work,
out of these measures, Kharal had stated that only the Euclidean distance
and the Normalized Euclidean distance between two T1SS were metrics
(refer to $Definition\:19,Proposition\:20$, pages $8-9$)\cite{kdist}.
These two distances were respectively defined in\cite{kdist} as\\
\\
Euclidean distance\cite{kdist}:\\
\\
 $e((F,A),(G,B))=|A\triangle B|+\sqrt{\sum_{\alpha\epsilon A\cap B}|F(\alpha)\triangle G(\alpha)|^{2}}$\\
\\
Normalized Euclidean distance\cite{kdist}:\\
\\
$q((F,A),(G,B))=\frac{|A\triangle B|}{\sqrt{|A\cup B|}}+\sqrt{\sum_{\alpha\epsilon A\cap B}\chi(\alpha)}$\\
where $\chi(\alpha)=\begin{cases}
\frac{|F(\alpha)\triangle G(\alpha)|^{2}}{|F(\alpha)\cup G(\alpha)|} & if\, F(\alpha)\cup G(\alpha)\neq\varphi\\
0 & otherwise
\end{cases}$\\

However, in this respect, a deeper study reveals that these two above
mentioned measures of distance contain fallacies since they do not
always satisfy the triangle inequality $(M3)$. In order to establish
our point, we provide a counter-example in support of our argument
as follows:\\
\\
\textbf{Example 3.2.} Let $(F,A)$, $(G,B)$ and$(H,C)$ be T1SS defined
over the universe $U=\{x_{1},x_{2},x_{3},x_{4},x_{5}\}$ and the set
of parameters $E=\{\alpha_{1},\alpha_{2},\alpha_{3}\}$ such that:\\
$(F,A)=\{(\alpha_{1},\{x_{1},x_{2}\})\}$;$(G,B)=\{(\alpha_{2},\{x_{2},x_{3}\}),(\alpha_{3},\{x_{1},x_{4}\})\}$;\\
$(H,C)=\{(\alpha_{1},\{x_{3},x_{4},x_{5}\}),(\alpha_{2},\{x_{2},x_{3}\}),(\alpha_{3},\{x_{1},x_{3},x_{4}\})\}$\\

On calculation we get $e((F,A),(G,B))=3$, $e((G,B),(H,C))=2$ and
$e((F,A),(H,C))=7$. i.e. $e((F,A),(H,C))>e((F,A),(G,B))+e((G,B),(H,C))$.
Also, $q((F,A),(G,B))=1.155$, $e((G,B),(H,C))=1.155$ and $e((F,A),(H,C))=3.391$
i.e. $e((F,A),(H,C))>e((F,A),(G,B))+e((G,B),(H,C))$.

Thus, the Euclidean and the Normalized Euclidean distances, which
were referred to as metrics are, in reality, not metrics since they
violate the triangle inequality. \\
\\
In view of the above situation, we hereby propose some distance measures
between two T1SS as\\
\\
\textbf{Definition 3.3.} The parameter based distance measure is defined
as\\
$d_{p}((F,A),(G,B))=|A\cup B|-|A\cap B|+|F^{\#}\cup G^{\#}|-|F^{\#}\cap G^{\#}|$
where $F^{\#}=\cup_{\alpha\epsilon A}\{F(\alpha)\}$ and $G^{\#}=\cup_{\beta\epsilon B}\{G(\beta)\}$\\
\\
\textbf{Theorem 3.4.} $d_{p}((F,A),(G,B))$ is a metric.\\
\\
Proof: We only give the proof of the Triangle inequality since the
other proofs are straight-forward.\\
Consider any three T1SS $(F,A),(G,B),(H,C)\epsilon S(X,E).$Now,\\
$\left|A\cup B\right|-\left|A\cap B\right|+\left|B\cup C\right|-\left|B\cap C\right|-\left|A\cup C\right|+\left|A\cap C\right|$\\
$=2(|B|-|A\cap B|-|B\cap C|+|A\cap C|)$\\
Also, $B=(B\cap(A\cup C))\cup(B\cap(A\cup C)^{c})$\\
$\Rightarrow|B|=|B\cap(A\cup C)|+|B\cap(A\cup C)^{c}|=|(A\cap B)\cup(B\cap C)|+|B\cap(A^{c}\cap C^{c})|=|A\cap B|+|B\cap C|-|A\cap B\cap C|+|B\cap A^{c}\cap C^{c}|$\\
i.e. $|B|-|A\cap B|-|B\cap C|=-|A\cap B\cap C|+|B\cap A^{c}\cap C^{c}|$\\
So, $2(|B|-|A\cap B|-|B\cap C|+|A\cap C|)=2(|A\cap C|-|A\cap B\cap C|+|B\cap A^{c}\cap C^{c}|)\geq0$,
since $(A\cap B\cap C)\subseteq A\cap C$ \\
\\
Hence it follows that\\
 $|A\cup C|-|A\cap C|\leq|A\cup B|-|A\cap B|+|B\cup C|-|B\cap C|$\\
\\
Similarly, it can be proved that, \\
$|F^{\#}\cup H^{\#}|-|F^{\#}\cap H^{\#}|\leq|F^{\#}\cup G^{\#}|-|F^{\#}\cap G^{\#}|+|G^{\#}\cup H^{\#}|-|G^{\#}\cap H^{\#}|$\\
\\
\textbf{Definition 3.5.} The matrix-representation based distance
measure is defined as\\
$d_{m}((F,A),(G,B))=|A\cup B|-|A\cap B|+\sum_{\alpha\epsilon A\cup B}\sum_{x\epsilon X}|F(\alpha)(x)-G(\alpha)(x)|$,
where $F(\alpha)(x)=\begin{cases}
1 & if\: x\epsilon F(\alpha)\\
0 & otherwise
\end{cases}$; $G(\alpha)(x)=\begin{cases}
1 & if\: x\epsilon G(\alpha)\\
0 & otherwise
\end{cases}$\\
\\
\textbf{Theorem 3.6.} The measure $d_{m}$ is a metric.\\
\\
Proof: We only give an outline of the proof of the triangle inequality
since the rest of the proofs are straight-forward:\\
\\
Consider any three T1SS $(F,A),(G,B),(H,C)\epsilon S(X,E).$ Proceeding
in an exactly same way as the proof of Theorem 3.4. , we can show
that for any three parameter sets $A,B,C$ \\
 $|A\cup C|-|A\cap C|\leq|A\cup B|-|A\cap B|+|B\cup C|-|B\cap C|$\\
\\
Also, $|F(\alpha)(x_{i})-H(\alpha)(x_{i})|=|F(\alpha)(x_{i})-G(\alpha)(x_{i})+G(\alpha)(x_{i})-H(\alpha)(x_{i})|\leq|F(\alpha)(x_{i})-G(\alpha)(x_{i})|+|G(\alpha)(x_{i})-H(\alpha)(x_{i})|$\\
$\Rightarrow d_{m}((F,A),(H,C))\leq d_{m}((F,A),(G,B))+d_{m}((G,B),(H,C))$\\
\\
We now proceed to define the distance measure between two T2SS as
a generalization of the distance measures for T1SS.\\
\\
\textbf{Definition} \textbf{3.7.} A mapping $D:S_{2}(X,E)\times S_{2}(X,E)\rightarrow\mathbb{R^{\dotplus}}$,
where $S_{2}(X,E)$ denotes the set of all type-2 soft sets over the
soft universe $(X,E)$, is said to be a distance measure between any
two T2SS if and only if for all $[\mathscr{F},A],[\mathscr{G},B],[\mathscr{H},C]\epsilon S_{2}(X,E)$
it satisfies the following conditions:\\
$(d1)$ $D([\mathscr{F},A],[\mathscr{G},B])=D([\mathscr{G},B],[\mathscr{F},A])$
\\
$(d2)$ $D([\mathscr{F},A],[\mathscr{G},B])\geq0$\\
$(d3)$ $D([\mathscr{F},A],[\mathscr{G},B])=0$ iff $[\mathscr{F},A]=[\mathscr{G},B]$.\\
$(d4)$ $D([\mathscr{F},A],[\mathscr{H},C])\leq D([\mathscr{F},A],[\mathscr{G},B])+D([\mathscr{G},B],[\mathscr{H},C])$\\
In addition to the above conditions if a distance measure satisfies
the following property viz.\\
$(d5)$ $D([\mathscr{F},A],[\mathscr{G},B])\leq1$\\
it is said to be a Normalized distance measure.\\
\\
\textbf{Theorem 3.8.} The parameter based distance measure for T2SS
defined as\\
$D_{p}([\mathscr{F},A],[\mathscr{G},B])=|A\cup B|-|A\cap B|+|E_{A}\cup E_{B}|-|E_{A}\cap E_{B}|+|F^{\#\#}\cup G^{\#\#}|-|F^{\#\#}\cap G^{\#\#}|$\\
where $F^{\#\#}=\cup_{\alpha\epsilon A}\cup_{\beta\epsilon E_{A}^{\alpha}}\{F_{\alpha}(\beta)\}$
and $G^{\#\#}=\cup_{\alpha\epsilon B}\cup_{\beta\epsilon E_{B}^{\alpha}}\{G_{\alpha}(\beta)\}$\\
is a metric.\\
\\
Proofs are similar to those of Theorem 3.4.\\
\\
\textbf{Theorem} \textbf{3.9.} The mapping $D_{m}:S_{2}(X,E)\times S_{2}(X,E)\rightarrow\mathbb{R^{\dotplus}}$
defined by\\
\\
$D_{m}([\mathscr{F},A],[\mathscr{G},B])=|A\cup B|-|A\cap B|+|E_{A}\cup E_{B}|-|E_{A}\cap E_{B}|$\\
$+\sum_{\alpha\epsilon A\cap B}\sum_{\beta\epsilon E_{A}^{\alpha}\cap E_{B}^{\alpha}}\sum_{x\epsilon X}|F_{\alpha}(\beta)(x)-G_{\alpha}(\beta)(x)|$\\
where $F_{\alpha}(\beta)(x)=\begin{cases}
1 & if\: x\epsilon F_{\alpha}(\beta)\\
0 & otherwise
\end{cases}$; $G_{\alpha}(\beta)(x)=\begin{cases}
1 & if\: x\epsilon G_{\alpha}(\beta)\\
0 & otherwise
\end{cases}$\\
is a distance measure between the T2SS $[\mathscr{F},A]$ and $[\mathscr{G},B]$.
It is the matrix-representation based distance measure between T2SS.\\
\\
\textbf{Theorem} \textbf{3.10.} The mappings $ND_{p}:S_{2}(X,E)\times S_{2}(X,E)\rightarrow[0,1]$
and $ND_{m}:S_{2}(X,E)\times S_{2}(X,E)\rightarrow[0,1]$ defined
as\\
\\
$ND_{p}([\mathscr{F},A],[\mathscr{G},B])=\frac{1}{|X|\times|A\cup B|\times|E_{A}\cup E_{B}|}\times D_{p}([\mathscr{F},A],[\mathscr{G},B])$\\
\\
$ND_{m}([\mathscr{F},A],[\mathscr{G},B])=\frac{1}{|X|\times|A\cup B|\times|E_{A}\cup E_{B}|}\times D_{m}([\mathscr{F},A],[\mathscr{G},B])$\\
\\
are distance measures. Furthermore, these measures are the normalized
parameter-based and normalized matrix-representation based distance
measures for T2SS.\\
\\
\textbf{Remark 3.11.} If, in particular, for $[\mathscr{F},A],[\mathscr{G},B],[\mathscr{H},C]\epsilon S_{2}(X,E)$,
$[\mathscr{F},A]\sqsubseteq[\mathscr{G},B]\sqsubseteq[\mathscr{H},C]$,
then the following results hold\\
$(i)$ $D_{p}([\mathscr{F},A],[\mathscr{H},C])=D_{p}([\mathscr{F},A],[\mathscr{G},B])+D_{p}([\mathscr{G},B],[\mathscr{H},C])$\\
$(ii)$ $D_{m}([\mathscr{F},A],[\mathscr{H},C])=D_{m}([\mathscr{F},A],[\mathscr{G},B])+D_{m}([\mathscr{G},B],[\mathscr{H},C])$\\
$(iii)$ $ND_{p}([\mathscr{F},A],[\mathscr{H},C])=ND_{p}([\mathscr{F},A],[\mathscr{G},B])+ND_{p}([\mathscr{G},B],[\mathscr{H},C])$\\
$(iv)$ $ND_{m}([\mathscr{F},A],[\mathscr{H},C])=ND_{m}([\mathscr{F},A],[\mathscr{G},B])+ND_{m}([\mathscr{G},B],[\mathscr{H},C])$\\
\textbf{}\\
\textbf{Example 3.12.} Consider the T2SS stated in example 2.7. Then,\\
$D_{p}([\mathscr{F},A],[\mathscr{G},B])=6$ units, $D_{m}([\mathscr{F},A],[\mathscr{G},B])=0.08$
units, $ND_{p}([\mathscr{F},A],[\mathscr{G},B])=7$ units and $Nd_{m}([\mathscr{F},A],[\mathscr{G},B])=0.093$
units.

\section{Entropy measure}

In this section we propose the definition of entropy and an entropy
measure for T2SS. \\
\textbf{}\\
\textbf{Definition} \textbf{4.1.} A T2SS $[\mathscr{F},A]$ is said
to be equivalent to another T2SS $[\mathscr{G},B]$ if there exists
two bijective functions $\psi:A\rightarrow B$ and $\Gamma:E_{A}\rightarrow E_{B}$
such that $\psi(\mathscr{F}_{x}^{*})=\mathscr{G}_{x}^{*}$ and $\Gamma(\mathscr{F}_{x}^{**})=\mathscr{G}_{x}^{**}$,
for all $x\epsilon X$, where\\
$\mathscr{F}_{x}^{*}=\{\alpha\epsilon A:x\epsilon F_{\alpha}(\beta)\: for\, some\,\beta\epsilon E_{A}^{\alpha}\}$,
$\mathscr{F}_{x}^{**}=\cup_{\alpha\epsilon A}\{\beta\epsilon E_{A}^{\alpha}:x\epsilon F_{\alpha}(\beta)\}$,
\\
$\mathscr{G}_{x}^{*}=\{\alpha\epsilon B:x\epsilon G_{\alpha}(\beta)\: for\, some\,\beta\epsilon E_{B}^{\alpha}\}$,
$\mathscr{G}_{x}^{**}=\cup_{\alpha\epsilon B}\{\beta\epsilon E_{B}^{\alpha}:x\epsilon G_{\alpha}(\beta)\}$.\\
If a T2SS $[\mathscr{G},B]$ is equivalent to $[\mathscr{F},A]$,
it is written as, $[\mathscr{G},B]\backsimeq[\mathscr{F},A]$. \\
\\
\textbf{Remark 4.2.} In particular, if for any two T2SS $[\mathscr{F},A],[\mathscr{G},B]$,
$[\mathscr{F},A]\backsimeq[\mathscr{G},B]$ then $|A|=|B|$, $|E_{A}|=|E_{B}|$,
$|\mathscr{F}_{x}^{*}|=|\mathscr{G}_{x}^{*}|$ and $|F_{x}^{**}|=|\mathscr{G}_{x}^{**}|$.\\
\\
\textbf{Definition} \textbf{4.3.} A T2SS $[\mathscr{F},A]$ is said
to be deterministic if\\
$(i)$ $A\cap E_{A}=\varphi$\\
$(ii)$ for distinct $\alpha_{1},\alpha_{2}\epsilon A$ (i.e. $\alpha_{1}\neq\alpha_{2}$),
$E_{A}^{\alpha_{1}}\cap E_{A}^{\alpha_{2}}=\varphi$ and$\mathscr{F}(\alpha_{1})\tilde{\cap}\mathscr{F}(\alpha_{2})=\tilde{\varphi}$,
the null type-1 soft set.\\
$(iii)$ for distinct $\beta_{1},\beta_{2}\epsilon E_{A}^{\alpha}$
(i.e. $\beta_{1}\neq\beta_{2}$), $F_{\alpha}(\beta_{1})\cap F_{\alpha}(\beta_{2})=\varphi$.\\
$(iv)$ $\cup_{\alpha\epsilon A}\{\cup_{\beta\epsilon E_{A}^{\alpha}}F_{\alpha}(\beta)\}=X$\\
\\
\textbf{Example 4.4.} Consider a T2SS $[\mathscr{F},A]$, over the
crisp universe $X=\left\{ x_{1},x_{2},x_{3},x_{4},x_{5}\right\} $
and the sets of parameters $A=\left\{ \alpha_{1},\alpha_{2},\alpha_{3}\right\} $,
$E_{A}=\left\{ \beta_{1},\beta_{2},\beta_{3},\beta_{4}\right\} $
defined as,\\
$\mathscr{F}\left(\alpha_{1}\right)=\left\{ \frac{\beta_{1}}{\left\{ x_{1},x_{2}\right\} }\right\} $;$\mathscr{F}\left(\alpha_{2}\right)=\left\{ \frac{\beta_{2}}{\left\{ x_{4}\right\} },\frac{\beta_{4}}{\left\{ x_{5}\right\} }\right\} $;$\mathscr{F}\left(\alpha_{3}\right)=\left\{ \frac{\beta_{3}}{\left\{ x_{3}\right\} }\right\} $\\
The T2SS $\left[\mathscr{F},A\right]$, mentioned above is a deterministic
T2SS.\\
\\
\textbf{Remark 4.5.} The null and the absolute T2SS, respectively
denoted by$\tilde{\Phi}$and $\mathbb{\tilde{A}}$ are not deterministic
T2SS.\\
 \\
\textbf{Definition} \textbf{4.6.} A mapping $E_{m}:S_{2}(X,E)\rightarrow\mathbb{R^{\dotplus}}$,
where $S_{2}(X,E)$ denotes the set of all type-2 soft sets over the
soft universe $(X,E)$, is said to be a measure of entropy of a T2SS
if and only if it satisfies the following conditions:\\
$(e1)$ $E_{m}(\tilde{\Phi})=1$, $E_{m}(\tilde{\mathbb{A}})=1$\\
$(e2)$ $E_{m}([\mathscr{F},A])\leq E_{m}([\mathscr{G},B])$ if $[\mathscr{F},A](\neq\tilde{\Phi})\sqsubseteq[\mathscr{G},B]$.\\
$(e3)$ $E_{m}([\mathscr{F},A])=0$ if $[\mathscr{F},A]$ is a deterministic
T2SS. \\
$(e4)$ $E_{m}([\mathscr{F},A])=E_{m}([\mathscr{G},B])$ iff $[\mathscr{F},A]\backsimeq[\mathscr{G},B]$.\\
\\
\textbf{Definition} \textbf{4.7.} For a T2SS $[\mathscr{F},A]$, define
$E_{m}:S_{2}(X,E)\rightarrow\mathbb{R^{\dotplus}}$ as\\
$E_{m}([\mathscr{F},A])=\begin{cases}
1 & if\:[\mathscr{F},A]=\tilde{\mathbb{A}\:}or\:[\mathscr{F},A]=\tilde{\Phi}\\
1-\frac{2|X|}{\sum_{x\epsilon X}\{|\mathscr{F}_{x}^{*}|+|\mathscr{F}_{x}^{**}|\}} & otherwise
\end{cases}$\\
\textbf{}\\
\textbf{Theorem} \textbf{4.8.} $E_{m}([\mathscr{F},A])$ as defined
above, is a measure of entropy for T2SS.\\
\\
\textbf{Example 4.9.} For example 2.7, the entropy measure $E_{m}([\mathscr{F},A])=0.44$.

\section{Similarity measures}

\textbf{Definition} \textbf{5.2.} For two T2SS $[\mathscr{F},A],[\mathscr{G},B]$,
define the soft real number valued similarity measure as a soft real
valued mapping $S^{'}:(A\cup B)\rightarrow\mathscr{P}(R)$, over the
set of parameters $A\cup B$, such that for any two T2SS $[\mathscr{F},A],[\mathscr{G},B]\epsilon S_{2}(X,E)$
that are not disjoint, for each $\alpha\epsilon A\cup B$ ,$S^{'}(\alpha)$
is defined as\\
$S^{'}(\alpha)=S_{\alpha}((F_{\alpha},E_{A}^{\alpha}),(G_{\alpha},E_{B}^{\alpha}))$
where $S_{\alpha}((F_{\alpha},E_{A}^{\alpha}),(G_{\alpha},E_{B}^{\alpha}))$
is the similarity measure between the T1SS corresponding to each primary
parameter $\alpha$ i.e. it satisfies the conditions:\\
$(s^{'}1)$ $S_{\alpha}((F_{\alpha},E_{A}^{\alpha}),(G_{\alpha},E_{B}^{\alpha}))=S_{\alpha}((G_{\alpha},E_{B}^{\alpha}),(F_{\alpha},E_{A}^{\alpha}))$\\
$(s^{'}2)$ $0\leq S_{\alpha}((F_{\alpha},E_{A}^{\alpha}),(G_{\alpha},E_{B}^{\alpha}))\leq1$\\
$(s^{'}3)$ $S_{\alpha}((F_{\alpha},E_{A}^{\alpha}),(G_{\alpha},E_{B}^{\alpha}))=1$
iff $(F_{\alpha},E_{A}^{\alpha})=(G_{\alpha},E_{B}^{\alpha})$\\
$(s^{'}4)$ for$(F_{\alpha},E_{A}^{\alpha}),(G_{\alpha},E_{B}^{\alpha}),(H_{\alpha},E_{C}^{\alpha})\epsilon S(X,E)$,
where $S(X,E)$ denotes the collection of all T1SS over the initial
soft universe $(X,E)$, if $(F_{\alpha},E_{A}^{\alpha})\tilde{\subseteq}(G_{\alpha},E_{B}^{\alpha})\tilde{\subseteq}(H_{\alpha},E_{C}^{\alpha})$
then $S_{\alpha}((F_{\alpha},E_{A}^{\alpha}),(H_{\alpha},E_{C}^{\alpha}))\leq S_{\alpha}((F_{\alpha},E_{A}^{\alpha}),(G_{\alpha},E_{B}^{\alpha}))\wedge S_{\alpha}((G_{\alpha},E_{B}^{\alpha}),(H_{\alpha},E_{C}^{\alpha}))$\\
$(s^{'}5)$ $S_{\alpha}((F_{\alpha},E_{A}^{\alpha}),(G_{\alpha},E_{B}^{\alpha}))=0$,
for all $\alpha\epsilon A\cup B$ if $[\mathscr{F},A],[\mathscr{G},B]\epsilon S_{2}(X,E)$
are disjoint.\\
\\
\textbf{Definition} \textbf{5.3.} Suppose $[\mathscr{F},A],[\mathscr{G},B]\epsilon S_{2}(X,E)$.
Define for $\alpha\epsilon A\cup B$\\
$S(\alpha)=\begin{cases}
\frac{1}{|E_{A}^{\alpha}\cup E_{B}^{\alpha}|}\times\left\{ \sum_{\beta\epsilon E_{A}^{\alpha}\cap E_{B}^{\alpha}}\frac{|F_{\alpha}(\beta)\cap G_{\alpha}(\beta)|}{|F_{\alpha}(\beta)\cup G_{\alpha}(\beta)|}\right\}  & when\,\alpha\epsilon A\cap B\neq\phi,E_{A}^{\alpha}\cap E_{B}^{\alpha}\neq\varphi\\
0 & otherwise
\end{cases}$\\
\\
\textbf{Theorem} \textbf{5.4.} $S$ is a soft real number valued similarity
measure between $[\mathscr{F},A]$ and $[\mathscr{G},B]$.\\
\\
\textbf{Definition 5.5.} A mapping $S_{m}:S_{2}(X,E)\times S_{2}(X,E)\rightarrow\mathbb{R^{\dotplus}}$,
where $\mathbb{R^{\dotplus}}$ denotes the set of all positive reals
and $S_{2}(X,E)$ denotes the set of all type-2 soft sets over the
soft universe $(X,E)$, is said to be a measure of similarity between
any two T2SS if and only if for all $[\mathscr{F},A],[\mathscr{G},B],[\mathscr{H},C]\epsilon S_{2}(X,E)$
it satisfies the following conditions:\\
$(s1)$ $S_{m}([\mathscr{F},A],[\mathscr{G},B])=S_{m}([\mathscr{G},B],[\mathscr{F},A])$\\
$(s2)$ $0\leq S_{m}([\mathscr{F},A],[\mathscr{G},B])\leq1$ \\
$(s3)$ $S_{m}([\mathscr{F},A],[\mathscr{G},B])=1$ iff $[\mathscr{F},A]=[\mathscr{G},B]$.\\
$(s4)$ if $[\mathscr{F},A]\sqsubseteq[\mathscr{G},B]\sqsubseteq[\mathscr{H},C]$
then $S_{m}([\mathscr{F},A],[\mathscr{H},C])\leq S_{m}([\mathscr{F},A],[\mathscr{G},B])\land S_{m}([\mathscr{G},B],[\mathscr{H},C])$\\
\\
\textbf{Theorem} \textbf{5.6.} The mapping $S_{m}:S_{2}(X,E)\times S_{2}(X,E)\rightarrow\mathbb{R^{\dotplus}}$
defined as\\
$S_{m}([\mathscr{F},A],[\mathscr{G},B])=\begin{cases}
\frac{1}{|A\cup B|}\times\left\{ \sum_{\alpha\epsilon A\cap B}\{S(\alpha)\}\right\}  & when\,\alpha\epsilon A\cap B\neq\phi\\
0 & otherwise
\end{cases}$ \\
is a real valued similarity measure between the T2SS $[\mathscr{F},A],[\mathscr{G},B]\epsilon S_{2}(X,E)$,
where $S(\alpha)$ is a soft real valued similarity measure between
$[\mathscr{F},A],[\mathscr{G},B]$.\\
\\
Proof: We give the proof of $(s4)$. Other proofs can be done similarly.\\
\\
$(s4)$ Suppose, for $[\mathscr{F},A],[\mathscr{G},B],[\mathscr{H},C]\epsilon S_{2}(X,E)$,
such that, $[\mathscr{F},A]\sqsubseteq[\mathscr{G},B]\sqsubseteq[\mathscr{H},C]$.
Thus, from the conditions of containment for T2SS,\\
$A\subseteq B\subseteq C$, $\mathscr{F}(\alpha)\tilde{\subset}\mathscr{G}(\alpha)$,
for each $\alpha\epsilon A$ and $\mathscr{G}(\beta)\tilde{\subset}\mathscr{H}(\beta)$,
for each $\beta\epsilon B$.\\
Consider, $S_{m}([\mathscr{F},A],[\mathscr{H},C])$ and $S_{m}([\mathscr{F},A],[\mathscr{G},B])$.\\
$\frac{|A\cap C|}{|A\cup C|}=\frac{|A|}{|C|}$ and $\frac{|A\cap B|}{|A\cup B|}=\frac{|A|}{|B|}$.
Since, $B\subseteq C$, $|B|\leq|C|$. Thus, \\
$\frac{|A|}{|C|}\leq\frac{|A|}{|B|}$ i.e. $\frac{|A\cap C|}{|A\cup C|}\leq\frac{|A\cap B|}{|A\cup B|}$.\\
Again, since $\mathscr{F}(\alpha)\tilde{\subset}\mathscr{G}(\alpha)\tilde{\subset}\mathscr{H}(\alpha)$,
for each $\alpha\epsilon A$,$E_{A}^{\alpha}\subseteq E_{B}^{\alpha}\subseteq E_{C}^{\alpha}$.
Hence, $\frac{|E_{A}^{\alpha}\cap E_{C}^{\alpha}|}{|E_{A}^{\alpha}\cup E_{C}^{\alpha}|}=\frac{|E_{A}^{\alpha}|}{|E_{C}^{\alpha}|}\leq\frac{|E_{A}^{\alpha}\cap E_{B}^{\alpha}|}{|E_{A}^{\alpha}\cup E_{B}^{\alpha}|}=\frac{|E_{A}^{\alpha}|}{|E_{B}^{\alpha}|}$.\\
Also, $|E_{A}^{\alpha}\cap E_{C}^{\alpha}|=|E_{A}^{\alpha}|=|E_{A}^{\alpha}\cap E_{B}^{\alpha}|$
$\Rightarrow\sum_{\alpha\epsilon A\cap C}|E_{A}^{\alpha}\cap E_{C}^{\alpha}|=\sum_{\alpha\epsilon A\cap B}|E_{A}^{\alpha}\cap E_{B}^{\alpha}|$\\
Now, $F_{\alpha}(\beta)\subseteq G_{\alpha}(\beta)\subseteq H_{\alpha}(\beta)$,
$\alpha\epsilon A\cap B\cap C$. Thus, $\frac{|F_{\alpha}(\beta)\cap H_{\alpha}(\beta)|}{|F_{\alpha}(\beta)\cup H_{\alpha}(\beta)|}=\frac{|F_{\alpha}(\beta)|}{|H_{\alpha}(\beta)|}\leq\frac{|F_{\alpha}(\beta)|}{|G_{\alpha}(\beta)|}=\frac{|F_{\alpha}(\beta)\cap G_{\alpha}(\beta)|}{|F_{\alpha}(\beta)\cup G_{\alpha}(\beta)|}$.
Hence,\\
$\frac{|A\cap C|}{|A\cup C|}\times\frac{\sum_{\alpha\epsilon A\cap C}\{\frac{|E_{A}^{\alpha}\cap E_{C}^{\alpha}|}{|E_{A}^{\alpha}\cup E_{C}^{\alpha}|}\times\frac{\sum_{\beta\epsilon E_{A}\cap E_{C}}\frac{|F_{\alpha}(\beta)\cap H_{\alpha}(\beta)|}{|F_{\alpha}(\beta)\cup H_{\alpha}(\beta)|}}{\sum_{\alpha\epsilon A\cap C}|E_{A}^{\alpha}\cap E_{C}^{\alpha}|}\}}{|A\cap C|}$\\
$=\frac{|A|}{|C|}\times\frac{\sum_{\alpha\epsilon A}\{\frac{|E_{A}^{\alpha}|}{|E_{C}^{\alpha}|}\times\frac{\sum_{\beta\epsilon E_{A}}\frac{|F_{\alpha}(\beta)|}{|H_{\alpha}(\beta)|}}{\sum_{\alpha\epsilon A}|E_{A}^{\alpha}|}\}}{|A|}$\\
$\leq\frac{|A|}{|B|}\times\frac{\sum_{\alpha\epsilon A}\{\frac{|E_{A}^{\alpha}|}{|E_{B}^{\alpha}|}\times\frac{\sum_{\beta\epsilon E_{A}}\frac{|F_{\alpha}(\beta)|}{|G_{\alpha}(\beta)|}}{\sum_{\alpha\epsilon A}|E_{A}^{\alpha}|}\}}{|A|}$\\
$=\frac{|A\cap B|}{|A\cup B|}\times\frac{\sum_{\alpha\epsilon A\cap C}\{\frac{|E_{A}^{\alpha}\cap E_{C}^{\alpha}|}{|E_{A}^{\alpha}\cup E_{C}^{\alpha}|}\times\frac{\sum_{\beta\epsilon E_{A}\cap E_{C}}\frac{|F_{\alpha}(\beta)\cap G\alpha(\beta)|}{|F_{\alpha}(\beta)\cup G_{\alpha}(\beta)|}}{\sum_{\alpha\epsilon A\cap C}|E_{A}^{\alpha}\cap E_{C}^{\alpha}|}\}}{|A\cap B|}$\\
$\Rightarrow S_{m}([\mathscr{F},A],[\mathscr{H},C])\leq S_{m}([\mathscr{F},A],[\mathscr{G},B])$.\\
In an exactly analogous manner, it can be shown that $S_{m}([\mathscr{F},A],[\mathscr{H},C])\leq S_{m}([\mathscr{G},B],[\mathscr{H},C])$
which completes the proof.

\subsection{Relation among similarity, distance and entropy measures for Type-2
soft sets}

\textbf{Theorem} \textbf{5.1.1.} If $D([\mathscr{F},A],[\mathscr{G},B])$
denotes any of the afore-mentioned proposed distance measure between
any two T2SS $[\mathscr{F},A]$ and $[\mathscr{G},B]$ then $S_{d}([\mathscr{F},A],[\mathscr{G},B])=\frac{1}{1+D([\mathscr{F},A],[\mathscr{G},B])}$
is a measure of similarity between the T2SS concerned.\\
\\
\textbf{Remark} \textbf{5.1.2.} The similarity measure $S_{d}([\mathscr{F},A],[\mathscr{G},B])$
derived in terms of the distance between two T2SS is termed as a distance-based
similarity measure for T2SS.\\
\\
\textbf{Theorem 5.1.3. }If $E_{m}([\mathscr{F},A])$ denotes any measure
of entropy for a T2SS $[\mathscr{F},A]$ then $S_{E}([\mathscr{F},A],[\mathscr{G},B])=1-|E_{m}([\mathscr{F},A]\sqcup[\mathscr{G},B])-E_{m}([\mathscr{F},A]\sqcap[\mathscr{G},B])|$
is a measure of similarity between the T2SS concerned.\textbf{}\\
\\
Only the proof of $(s3)$ is shown.\\
\\
$(s3)$ Since $0\leq E_{m}([\mathscr{F},A])\leq1$, for any $[\mathscr{F},A]\epsilon S_{2}(X,E)$,
it implies that $0\leq S_{E}([\mathscr{F},A],[\mathscr{G},B])\leq1$\\
Also, $1-|E_{m}([\mathscr{F},A]\sqcup[\mathscr{G},B])-E_{m}([\mathscr{F},A]\sqcap[\mathscr{G},B])|=0$
whenever $|E_{m}([\mathscr{F},A]\sqcup[\mathscr{G},B])-E_{m}([\mathscr{F},A]\sqcap[\mathscr{G},B])|=1$.
We have, from Remark 4.2,\\
$|A\cup B|=|A\cap B|$, $|(\mathscr{F}\sqcup\mathscr{G})_{x}^{*}|=|(\mathscr{F}\sqcap\mathscr{G})_{x}^{*}|$
and $|(\mathscr{F}\sqcup\mathscr{G})_{x}^{**}|=|(\mathscr{F}\sqcap\mathscr{G})_{x}^{**}|$.\\
 $|A\cup B|=|A\cap B|\Rightarrow A=B$...............................(1)\\
Again,\\
 $|(\mathscr{F}\sqcup\mathscr{G})_{x}^{*}|=|\{\alpha\epsilon A:x\epsilon F_{\alpha}(\beta_{1})\cup G_{\alpha}(\beta_{2}),\beta_{1},\beta_{2}\epsilon E_{A}^{\alpha}\cup E_{B}^{\alpha}\}|$\\
$=|(\mathscr{F}\sqcap\mathscr{G})_{x}^{*}|=|\{\alpha\epsilon A:x\epsilon F_{\alpha}(\beta_{1})\cap G_{\alpha}(\beta_{2}),\beta_{1},\beta_{2}\epsilon E_{A}^{\alpha}\cap E_{B}^{\alpha}\}|$.............(a)\\
 $|(\mathscr{F}\sqcup\mathscr{G})_{x}^{**}|=|\cup_{\alpha\epsilon A}\{\beta\epsilon E_{A}^{\alpha}:x\epsilon F_{\alpha}(\beta_{1})\cup G_{\alpha}(\beta_{2}),\beta_{1},\beta_{2}\epsilon E_{A}^{\alpha}\cup E_{B}^{\alpha}\}|$\\
$=|(\mathscr{F}\sqcap\mathscr{G})_{x}^{**}|=|\cup_{\alpha\epsilon A}\{\beta\epsilon E_{A}^{\alpha}:x\epsilon F_{\alpha}(\beta_{1})\cup G_{\alpha}(\beta_{2}),\beta_{1},\beta_{2}\epsilon E_{A}^{\alpha}\cap E_{B}^{\alpha}|$...........(b)\\
 From (b) we have, $|E_{A}^{\alpha}\cup E_{B}^{\alpha}|=|E_{A}^{\alpha}\cap E_{B}^{\alpha}|\Rightarrow E_{A}^{\alpha}=E_{B}^{\alpha}=E^{'},$say.....................(2)\\
From (a),(b) we can conclude, $|\{x:x\epsilon F_{\alpha}(\beta)\cup G_{\alpha}(\beta),\beta\epsilon E^{'}\}|=|x:x\epsilon F_{\alpha}(\beta)\cap G_{\alpha}(\beta),\beta\epsilon E^{'}|$,
for each $\alpha\epsilon A.$\\
$\Rightarrow|F_{\alpha}(\beta)\cup G_{\alpha}(\beta)|=|F_{\alpha}(\beta)\cap G_{\alpha}(\beta)|$,
for each $\alpha\epsilon A$ and $\beta\epsilon E^{'}$ .\\
$\Rightarrow F_{\alpha}(\beta)=G_{\alpha}(\beta)$, for each $\alpha\epsilon A$
and $\beta\epsilon E^{'}.$\\
 Thus, $S_{E}([\mathscr{F},A],[\mathscr{G},B])=1\Leftrightarrow[\mathscr{F},A]=[\mathscr{G},B]$.\\
\\
\textbf{Example 5.1.4} Consider the T2SS stated in example 2.7. Then,
$S_{m}([\mathscr{F},A],[\mathscr{G},B])=0.167$. Also, the distance
based simiarity measures corresponding to the distance measure $D_{m}$
and normalized distance measure $ND_{m}$ are $0.125$ and $0.915$
respectively (refer example 3.12.).

\section{Application of the proposed soft real number valued similarity measure}

As for application of soft real valued similarity measure, we have
considered the case of a person sufferring from diabetes. Thus, quite
naturally the person would be advised to consume meals that are low
in carbohydrates and rich in proteins and fibres. We assume that out
of homely atmosphere, the person encounters different choices of food
items on the platter out of which he has to choose the best suited
food items which would not cause any harm to his present health status.Algorithm
for selecting the best suited food items from the available choices
involves the following steps:\\
\\
\textbf{Step 1:} Suppose that we have $n$ number of menus from $n$
pantries under consideration such that the available food are categorized
and listed in the form of $n$ number of T2SS $[\mathscr{F}_{i},A_{i}],i=1,2,...,n$
and let the ideal preference of food items be listed as T2SS $[\mathscr{F},A]$.\\
\textbf{Step 2:} Calculate the soft real number valued similarity
measures $S^{'}([\mathscr{F},A],[\mathscr{F}_{i},A_{i}])$ for $i=1,2,...,n$.\\
\textbf{Step 3:} For a primary parameter $\alpha\epsilon A$, compare
the values of $S^{'}([\mathscr{F},A],[\mathscr{F}_{1},A_{1}])(\alpha)$,
$S^{'}([\mathscr{F},A],[\mathscr{F}_{2},A_{2}])(\alpha)$, $...$
and $S^{'}([\mathscr{F},A],[\mathscr{F}_{i},A_{i}])(\alpha)$ and
select\\
$max_{i}\{S^{'}([\mathscr{F},A],[\mathscr{F}_{i},A_{i}])(\alpha)\}$.
\\
If for $i=k\epsilon[1,n]$, $S^{'}([\mathscr{F},A],[\mathscr{F}_{k},A_{k}])(\alpha)=max_{i}\{S^{'}([\mathscr{F},A],[\mathscr{F}_{i},A_{i}])(\alpha)\}$,
then the selected pantry is the $k^{th}$ pantry. Proceed to Step
4.\\
\textbf{Step 4:} The selected food items are obtained as $\mathscr{F}(\alpha)\tilde{\cap}\mathscr{F}_{k}(\alpha)$.
\\
\textbf{Step 5:} Repeat Step 3 for all $\alpha\epsilon A$.\\
\textbf{Step 6:} Stop.\\
\\
We consider a set of two different menus from two different pantries
and we proceed with representing the available food items in the respective
menus in terms of two different T2SS as follows, the primary classifying
parameters being the types of food available for breakfast, lunch,
dinner and supper.

Suppose, $X$ denote the set of all available food items given by\\
$X=\{pastry,bagels,brown\, bread,mousse,noodles,rice,fruit\, juice,cereals,pasta,$\\
$vegetables,club\, sandwich,chicken,salad,soup,fish,pudding,milk,fruits,egg,$\\
$chapati,nuts\}$.

Let, the entire set of underlying parameters is given by,\\
 $\{fibre\, rich,protein\, rich,carb.\, rich,soft\, diet,liquid\, diet\}$
and the primary set of parameters be $A=\{breakfast,lunch,dinner,supper\}$.

Suppose the food available at Pantry1 can be classified into a T2SS
$[\mathscr{F}_{1},A]$ as follows:\\
\\
$\mathscr{F}_{1}(breakfast)=\{\frac{carb.\, rich}{\{pastry,bagels\}},\frac{fluid\, diet}{\{fruit\, juice\}},\frac{fibre\, rich}{\{cereals,fruits\}}\}$\\
$\mathscr{F}_{1}(lunch)=\{\frac{carb.\, rich}{\{rice,noodles,pasta\}},\frac{protein\, rich}{\{fish\}},\frac{fibre\, rich}{\{vegetables,salad\}}\}$\\
$\mathscr{F}_{1}(dinner)=\{\frac{protein\, rich}{\{chicken\}},\frac{soft\, diet}{\{soup\}},\frac{fibre\, rich}{\{salad\}}\}$\\
$\mathscr{F}_{1}(supper)=\{\frac{carb.\, rich}{\{club\, sandwich\}},\frac{soft\, diet}{\{pudding\}}\}$

In a similar way, the food items available in Pantry2 is classified
into a T2SS $[\mathscr{F}_{2},A]$ as follows:\\
$\mathscr{F}_{2}(breakfast)=\{\frac{carb.\, rich}{\{bagels\}},\frac{protein\, rich}{\{egg,chicken\}},\frac{fluid\, diet}{\{fruit\, juice,milk\}},\frac{fibre\, rich}{\{brown\, bread,cereals,fruits\}}\}$\\
$\mathscr{F}_{2}(lunch)=\{\frac{carb.\, rich}{\{rice,noodles\}},\frac{fibre\, rich}{\{vegetables\}},\frac{soft\, diet}{\{mousse\}}\}$\\
$\mathscr{F}_{2}(dinner)=\{\frac{carb.\, rich}{\{noodles\}},\frac{protein\, rich}{\{fish\}}\frac{fibre\, rich}{\{chapati,vegetables\}}\}$\\
$\mathscr{F}_{2}(supper)=\{\frac{protein\, rich}{\{chicken\}},\frac{fibre\, rich}{\{nuts,salad\}}\}$

Finally we represent the ideal food items that the person is allowed
to eat throughout the day for the various course of meals in terms
of a T2SS $[\mathscr{F},A]$ as,\\
$\mathscr{F}(breakfast)=\{\frac{fibre\, rich}{\{cereals,fruits,brown\, bread\}},\frac{fluid\, diet}{\{milk\}}\}$\\
$\mathscr{F}(lunch)=\{\frac{protein\, rich}{\{fish,chicken,egg\}},\frac{fibre\, rich}{\{vegetables,salad,chapati\}}\}$\\
$\mathscr{F}(dinner)=\{\frac{protein\, rich}{\{chicken,fish\}},\frac{soft\, diet}{\{soup\}},\frac{fibre\, rich}{\{salad,vegetables,chapati\}}\}$\\
$\mathscr{F}(supper)=\{\frac{soft\, diet}{\{soup\}},\frac{fibre\, rich}{\{salad\}}\}$

We now attempt to solve our problem at hand.\\
 Here, $i=2$ and $A_{1}=A_{2}=A$.\\
The soft real valued similarity measures are given as,\\
$S^{'}([\mathscr{F},A],[\mathscr{F}_{1},A])=\{(breakfast,0.444),(lunch,0.333),(dinner,0.611),(supper,0.000)\}$\\
$S^{'}([\mathscr{F},A],[\mathscr{F}_{2},A])=\{(breakfast,0.375),(lunch,0.083),(dinner,0.292),(supper,0.167)\}$

Thus, comparing the values parameterwise we see that for the primary
parameter $breakfast$, $max_{i=1,2}\{S^{'}([\mathscr{F},A],[\mathscr{F}_{i},A_{i}])\}=S^{'}([\mathscr{F},A],[\mathscr{F}_{2},A])\}$.
We thus conclude that the person would opt for his breakfast from
Pantry2. Also, from further calculations we see that the person opts
for his $lunch$ from Pantry1, $dinner$ from Pantry1 and $supper$
from Pantry2.

Also, for the choice of food items we show the calculations pertaining
to the parameter viz. $breakfast$ as,\\
$\mathscr{F}(breakfast)\tilde{\cap}\mathscr{F}_{2}(breakfast)=\{\frac{fibre\, rich}{\{cereals,fruits,brown\, bread\}},\frac{fluid\, diet}{\{milk\}}\}$\\
Thus, the choice of food items for breakfast include $cereals,fruits,brown\, bread,$\\
$milk$. Similarly, for $lunch$,$dinner$ and $supper$, the list
of foods is given by,\\
$\mathscr{F}(lunch)\tilde{\cap}\mathscr{F}_{1}(lunch)=\{\frac{protein\, rich}{fish},\frac{fibre\, rich}{vegetables,salad}\}$\\
$\mathscr{F}(dinner)\tilde{\cap}\mathscr{F}_{1}(dinner)=\{\frac{protein\, rich}{chicken},\frac{soft\, diet}{soup},\frac{fibre\, rich}{salad}\}$\\
$\mathscr{F}(supper)\tilde{\cap}\mathscr{F}_{2}(supper)=\{\frac{fibre\, rich}{salad}\}$\\
i.e. for lunch, the permissible food items include $fish,vegetables,salad$;
the prefereable food for dinner include $chicken,soup,salad$ and
$salad$ for supper.\\
\\
\\
\textbf{Funding:} The research of the first author is supported by
University JRF (Junior Research Fellowship).\\
The research of the third author is partially supported by the Special
Assistance Programme (SAP) of UGC, New Delhi, India {[}Grant no. F
510/3/DRS-III/(SAP-I){]}.

\end{document}